\newtheorem{prop}{Proposition}
\newtheorem{theorem}{Theorem}

\documentstyle[12pt]{article}

\begin{document}
\title{The moduli space of special Lagrangian submanifolds}
\author{N.J.HITCHIN\\[5pt]
\itshape Department of Pure Mathematics and Mathematical Statistics\\
\itshape University of Cambridge\\
\itshape 16 Mill Lane\\
\itshape Cambridge CB2 1SB\\
\itshape England\\[12pt]}
\maketitle
\bigskip
\section{Introduction}
In their quest for examples of minimal submanifolds, Harvey and Lawson in 1982
\cite{HL}
extended the well-known fact that a complex submanifold of a K\"ahler manifold
is minimal to the more general context of {\it calibrated} submanifolds. One
such class  is that of  special Lagrangian submanifolds of a Calabi-Yau manifold.
New   developments  in the study of these have raised the question as to
whether they should be accorded equal status with complex submanifolds.
  
 The developments stem from two sources. The first is the deformation theory of
 R.C. McLean \cite{Mac}. This shows that, given one compact special Lagrangian
 submanifold $L$, there is a local moduli space which is a manifold and whose
 tangent space at $L$ is canonically identified with the space of harmonic
 $1$-forms  on $L$. The ${\cal L}^2$ inner product on harmonic forms then gives
 the  moduli  space  a natural  Riemannian metric. The second input is from the
 paper  of Strominger, Yau and  Zaslow \cite{SYZ} which studies the moduli space 
of special Lagrangian tori in
 the context of mirror symmetry.
 
 This paper is in some sense a commentary on these two works, but it is
 provoked by the question: ``What is the natural geometrical structure on the
 moduli space of special Lagrangian submanifolds in a Calabi-Yau manifold?''
 We know that a moduli space of complex submanifolds (when unobstructed) is a
 complex  manifold. We shall show that the moduli space $M$ of special
 Lagrangian  submanifolds   has the local structure of a  Lagrangian
 submanifold, and we conjecture that it
 is ``special'' in an appropriate sense. 
 
 ``A Lagrangian submanifold of what?'' the reader may well ask. Recall that if
 $V$ is a finite-dimensional real vector space, then the natural pairing with
 its dual space $V^*$ defines a symplectic structure on $V\times V^*$. It also
 defines an indefinite metric. We shall show that there is a natural embedding
 of the local moduli space $M$ as a Lagrangian submanifold in the product 
 $H^1(L,  {\bf  R}  )\times  H^{n-1}(L, {\bf R})$   (where $n=\dim L$) of two
 dual vector spaces 
 and  that McLean's metric is  the  natural induced metric. 
   
 The symplectic manifold $V\times V^*$ can be thought of in two ways as a
 cotangent bundle: 
 as  either  $T^*V$  or  $T^*V^*$. Thus the Lagrangian submanifold $M$ is
 defined  locally  as  the  graph of the derivative of a  function $\phi:
 V\rightarrow {\bf R}$ or $\psi: V^*\rightarrow {\bf R}$. We show that this 
 symmetry (which is really the Legendre transform)  lies behind the
 viewpoint in \cite{SYZ}, where it is viewed as a manifestation of mirror
 symmetry.  This involves studying the structure of the moduli space of
 Lagrangian submanifolds together with flat line bundles. We show that there is
 a natural complex structure and K\"ahler metric on this space, and that this
 is a Calabi-Yau metric if the embedding of $M$ above is special.

 \section{Calabi-Yau manifolds}

A {\it Calabi-Yau manifold} is a  K\"ahler manifold of complex dimension
$n$ with a covariant constant holomorphic $n$-form. Equivalently it is a
Riemannian  manifold with holonomy contained in $SU(n)$.

It is convenient for our purposes to play down the role of the complex structure
in describing such manifolds and to emphasise instead the role of three
closed forms, satisfying certain algebraic identities (see \cite{Sal}).
We have the K\"ahler 2-form $\omega$ and the real and imaginary parts
$\Omega_1$ and $\Omega_2$ of the covariant constant $n$-form. These satisfy some  
identities:
\vskip .25cm
\noindent (i) $\omega$ is non-degenerate

\noindent (ii) $\Omega_1+i\Omega_2$ is locally decomposable and non-vanishing

\noindent (iii) $\Omega_1 \wedge \omega=\Omega_2 \wedge \omega =0$

\noindent (iv) $(\Omega_1+i \Omega_2)\wedge(\Omega_1-i \Omega_2)=\omega^n$
(resp. $i\omega^n$) if $n$ is even (resp. odd)

\noindent (v) $d\omega=0,\quad d\Omega_1=0,\quad d\Omega_2=0$
\vskip .25cm
\noindent These conditions (together with a positivity condition) we now show 
serve to characterize Calabi-Yau manifolds. Firstly if
$\Omega^c=\Omega_1+i\Omega_2=$ is locally decomposable as $\theta_1 \wedge
\theta_2 \wedge...\wedge \theta_n$, then  take the  subbundle $\Lambda$ of 
$T^*M\otimes {\bf C}$
spanned by $\theta_1,\dots,\theta_n$.  By (iv) and the fact that $\omega^n\ne
0$, we have 
$$\theta_1\wedge\dots\wedge\theta_n\wedge\bar\theta_1\wedge\dots\wedge 
\bar\theta_n\ne 0$$
and so $T^*M=\Lambda+\bar \Lambda$ and we have an almost-complex
structure. In this description a $1$-form $\theta$ is of type $(1,0)$ if and
only if 
$\Omega^c\wedge \theta=0$. Since from (v) $d\Omega_1 = d\Omega_2=0$
this means that $\Omega^c\wedge d\theta=0$. Writing
\begin{equation}
d\theta=\sum a_{ij}\theta_i\wedge\theta_j +\sum
b_{ij}\theta_i\wedge\bar \theta_j+\sum c_{ij}\bar\theta_i\wedge\bar
\theta_j
\label{form}
\end{equation}
we see that $c_{ij}=0$. Thus the ideal generated by $\Lambda$ is closed under
exterior differentiation, and by the Newlander-Nirenberg theorem the structure
is integrable. 

Similarly, applying the decomposition of 2-forms (\ref{form}) to $\omega$,
(iii) implies that the $(0,2)$ component vanishes, and since $\omega$ is real,
it is of type $(1,1)$. It is closed by (v), so if the hermitian
form so defined is positive definite, then we have a K\"ahler metric.

Since $\Omega^c$ is closed and of type $(n,0)$ it is a non-vanishing
holomorphic section $s$ of the canonical bundle. Relative to the trivialization
$s$, the hermitian connection has connection form  given by $\partial \log
(\Vert s \Vert^2)$. But  property (iv) implies that it has constant length, so the
connection form vanishes and 
$s=\Omega^c$ is covariant constant.

\section{Special Lagrangian submanifolds}

A submanifold $L$ of a symplectic manifold $X$ is  Lagrangian if $\omega$
restricts to zero on $L$ and $\dim X=2 \dim L$. A submanifold of a Calabi-Yau 
manifold is  {\it special Lagrangian} if in addition $\Omega=\Omega_1$ restricts
to zero on $L$. This condition involves only two out of the three forms, and in 
many respects
what we shall be doing is to treat them both --- the 2-form $\omega$ and the
$n$-form $\Omega$ --- on the same footing. 
\vskip .25cm
\noindent {\it Remarks:}

\noindent 1. We could relax the definition a little since
$\Omega^c$ is a {\it chosen} holomorphic $n$-form: any constant multiple of 
$\Omega^c$
would also be covariant  constant, so under some circumstances we may need to say 
that $L$ is special
 Lagrangian if, for some non-zero $c_1,c_2\in {\bf R}$,
 $c_1\Omega_1+c_2\Omega_2=0$.
 
 \noindent 2. On a special Lagrangian submanifold $L$, the $n$-form $\Omega_2$
 restricts to a non-vanishing  form, so in particular $L$ is always
 oriented.
\vskip .25cm

\vskip .25cm
Examples of special Lagrangian submanifolds are difficult to find, and so
far consist of  three types:
\begin{itemize}
\item
Complex Lagrangian submanifolds of hyperk\"ahler manifolds
\item
Fixed points of a real structure on a Calabi-Yau manifold
\item
Explicit examples for non-compact Calabi-Yau manifolds
\end{itemize}
The hyperk\"ahler examples arise easily. In this case we have $n=2k$  and 
three K\"ahler forms $\omega_1,\omega_2,\omega_3$ corresponding to the three
complex structures $I,J,K$ of the hyperk\"ahler manifold. With respect to the 
complex
structure $I$ the form 
$\omega^c=(\omega_2+i\omega_3)$ is a holomorphic symplectic form. If $L$ is a 
complex Lagrangian submanifold (i.e. $L$ is a complex submanifold and $\omega^c$
vanishes on $L$), then the real and imaginary parts of this,  $\omega_2$ and 
$\omega_3$, vanish on $L$. Thus $\omega=\omega_2$ vanishes
and if $k$ is odd (resp. even), the real (resp. imaginary) part of
$\Omega^c=(\omega_3+i\omega_1)^k$ vanishes.
 Using the complex structure
$J$ instead of $I$, we see that $L$ is special Lagrangian. For examples here,
we can take any holomorphic curve in a K3 surface $S$, or its symmetric product in 
the Hilbert
scheme $S^{[m]}$, which is hyperk\"ahler from \cite{Bea}.
\vskip .15cm
If $X$ is a Calabi-Yau manifold with a real structure --- an antiholomorphic
involution $\sigma$
--- for which $\sigma^*\omega=-\omega$ and $\sigma^*\Omega=-\Omega$, then 
 the fixed point set (the set of real points of $X$) is easily seen to be a
 special  Lagrangian  submanifold $L$.
 \vskip .15cm
 All Calabi-Yau metrics on compact manifolds are produced by the existence
 theorem of Yau. In the non-compact case, Stenzel \cite{Sten} has some concrete
 examples. In particular $T^*S^n$ (with the complex structure of an affine
 quadric) has a complete Calabi-Yau metric for which the zero section is
 special Lagrangian. When $n=2$ this is the hyperk\"ahler Eguchi-Hanson metric.

 \section{Deformations of special Lagrangian submanifolds}

R.C.McLean has studied deformations of special Lagrangian submanifolds. His
main result is 
\begin{theorem} \cite{Mac}  A normal vector field $V$ to a compact
special Lagrangian submanifold $L$ is the deformation vector field to a normal
deformation through special Lagrangian submanifolds if and only if the
corresponding 1-form $IV$ on $L$ is harmonic. There are no obstructions to
extending a first order deformation to an actual deformation and the tangent space 
to such
deformations can be identified through the cohomology class of the harmonic
form with $H^1(L,{\bf R})$.
\end{theorem}

Let us briefly see how the tangent space to the (local) moduli space $M$ is
identified with the space of harmonic 1-forms. Consider a 1-parameter family
$L_t$ of Lagrangian submanifolds as a smooth map $f:{\cal
L}\rightarrow X$ of the manifold ${\cal
L}=L\times U$ to $X$ where $U\subset {\bf R}$ is an interval and $f(L,t)=L_t$.
Since each $L_t$ is Lagrangian, $f^*\omega$ restricts to zero on each fibre of
$p:{\cal L}\rightarrow U$ 
so  we can find a 1-form $\tilde\theta$ on ${\cal L}$ such that
$$f^*\omega=dt\wedge \tilde\theta$$
The restriction $\theta$ of $\tilde\theta$ to each fibre $L\times \{t\}$  is 
independent of the choice of $\tilde \theta$, and since $d\omega=0$,
it follows that  $$d\theta=0$$ 
Similarly, since $L_t$ is {\it special} Lagrangian, the $n$-form $\Omega$
vanishes on each fibre, so that
$$f^*\Omega=dt\wedge \tilde\varphi$$
and since $d\Omega=0$ we have $d\varphi=0$. Using the induced metric on $L_t$
one can show that 
$$\varphi=\ast \theta$$
so that $\theta$ is the required harmonic form.
\vskip .25cm
A more invariant way of seeing this is to take a section of the normal bundle
of $L_t$, since this is what an infinitesimal variation canonically describes.
Take a representative vector field $V$ on $X$ and form the interior product
$\iota(V)\omega$. Since $\omega$ vanishes on $L_t$, the restriction of 
$\iota(V)\omega$ 
to $L_t$ is a 1-form which is independent of the choice of $V$. Now
$df(\partial/\partial t)$ is  naturally a section of the normal bundle of
$L_t\subset X
$ and $\theta$ is then the corresponding 1-form.
\vskip .25cm
Suppose now we take local coordinates $t_1,\dots, t_m$ on the moduli space $M$
of deformations of $L=L_0$. Here of course, from McLean, we know that
$m=b_1(L)=\dim H^1(L)$. For each tangent vector $\partial/\partial t_j$ we
define as above a corresponding closed 1-form $\theta_j$ on
$L_t$ for each $t\in M$: 
$$\iota(\partial/\partial t_j)\omega=\theta_j$$
(with a slight abuse of notation).

Let $A_1,\dots, A_m$ be a basis for $H_1(L,{\bf Z})$ (modulo torsion), then we
can evaluate the closed form $\theta_j$ on the homology class $A_i$ to obtain a
period matrix $\lambda_{ij}$ which is a function on the moduli space:
$$\lambda_{ij}=\int_{A_i}\theta_j$$
Since by McLean's theorem, the harmonic forms $\theta_j$ are linearly
independent, it follows that $\lambda_{ij}$ is invertible. We can now
be explicit about the identification of the tangent space to $M$ with the 
cohomology group $H^1(L,{\bf R})$.
Let $\alpha_1,\dots,\alpha_m\in H^1(L,{\bf Z})$ be the  basis dual to
 $A_1,\dots,A_m$.
It follows that
\begin{equation}
\partial/\partial
t_j\mapsto [\iota(\partial/\partial
t_j)\omega]=\sum \lambda_{ij}\alpha_i
\label{coord}
\end{equation}
identifies $T_tM$ with $H^1(L, {\bf R})$.
\vskip .25cm
\noindent We now investigate further properties of the period matrix $\lambda$.
\begin{prop} \label{alpha}The 1-forms $\xi_i=\sum \lambda_{ij}dt_j$ on $M$
are closed. \end{prop}

\noindent{\it Proof:} We represent the full local family of deformations by a
map $f:{\cal M}\rightarrow X$ where ${\cal M}\cong L\times M$ with projection
$p:{\cal M}\rightarrow M$.
Choose smoothly in each fibre of $p$ a circle representing $A_i$ to give an
$n+1$-manifold ${\cal M}_i\subseteq {\cal M}$ fibering over $M$. Define the
1-form $\xi$ on $M$ by $$\xi=p_{*}f^*\omega $$
The push-down map $p_*$ (integration over the fibres) takes closed forms to
closed forms, and since $d\omega=0$, $df^*\omega=0$ and so $d\xi=0$.

Now in local coordinates $\omega=\sum_j dt_j\wedge \tilde\theta_j$ and 
$\tilde\theta_j$ restricts to $\theta_j$ on each fibre.
Since $\theta_j$ is closed,
integration over the fibres of ${\cal M}_i$ is just evaluation on the  homology
class $A_i$. Thus $\xi_i=\xi$ and $\xi_i$ is closed.
 \vskip .5cm
 
 From this Proposition, we can find on $M$ local functions $u_1,\dots,u_m$,
 well-defined up to the addition of a constant, such that
 \begin{equation}
 du_i=\xi_i=\sum_j \lambda_{ij}dt_j
 \label{u}
 \end{equation}
 Since $\lambda_{ij}$ is invertible, $u_1,\dots,u_m$ are local coordinates on
 $M$. More invariantly, we have a coordinate chart
 \begin{equation}
 u:M\rightarrow H^1(L,{\bf R})
 \label{uu}
 \end{equation}
 defined by
 $u(t)=\sum_iu_i\alpha_i$ which is independent of the choice of basis, and is
 well-defined up to a translation.
 \vskip .5cm
 Clearly, we should follow our even-handed policy with respect to $\omega$ and
 $\Omega$ and enact the same procedure for $\Omega$. Thus, the basis
 $\alpha_1,\dots,\alpha_m$ defines a basis $B_1,\dots,B_m$ of $H_{n-1}(L,{\bf
 Z})$  and we form a period matrix $\mu_{ij}$:
 $$\mu_{ij}=\int_{B_i}\varphi_j$$
 In a similar fashion we find local coordinates $v_1,\dots, v_m$ on $M$
 such that
 \begin{equation}
 dv_i=\sum_j \mu_{ij}dt_j
 \label{v}
 \end{equation}
 and 
 an invariantly defined map
 \begin{equation}
 v:M\rightarrow H^{n-1}(L,{\bf R})
 \label{vv}
 \end{equation}
 given, using the basis $\beta_1,\dots,\beta_m$ of $H^{n-1}(L,{\bf R})$ dual to
 $B_1,\dots, B_m$ 
 by    
 $v(t)=\sum_i v_i \beta_i$.
 \vskip .5cm
 \noindent We obtain from $u$ and $v$ a map 
 $$F:M\rightarrow H^1(L,{\bf R})\times H^{n-1}(L,{\bf R})$$
 defined by $F(t)=(u(t),v(t))$.
 \vskip .25cm
 Let us see now how this fits in with the natural ${\cal L}^2$ metric on $M$.
 Note that since $L$ is oriented, $H^1(L)$ and $H^{n-1}(L)$ are canonically
 dual. For any vector space $V$ there is a natural indefinite symmetric form on
 $V\oplus V^*$ defined by
 $$B((v,\alpha),(v,\alpha))=\langle v, \alpha \rangle$$
 Thus $H^1(L)\times H^{n-1}(L)$ has a natural flat indefinite metric $G$.
 \begin{prop} The ${\cal L}^2$ metric $g$ on $M$ is 
 $F^*G$.  \end{prop}
 \noindent{\it Proof:} From (\ref{coord}), we have
 $$dF(\partial/\partial t_j)=(\sum_i \lambda_{ij}\alpha_i,\sum_i 
\mu_{ij}\beta_i)$$
 Thus
 \begin{equation}
 F^*G(\sum_j a_j \partial/\partial t_j, \sum_j a_j
 \partial/\partial  t_j)= \sum_{i,j,k,l} a_j a_k \lambda_{ij}\mu_{lk} \langle
 \alpha_i, \beta_l \rangle=
 \sum_{i,j,k,l} a_j a_k \lambda_{ij}\mu_{lk} \int_L 
 \alpha_i \wedge \beta_l
 \label{fg} 
 \end{equation}
  But
 $$\int_L (\sum_i a_i \theta_i)\wedge \ast (\sum_i a_i \theta_i)= \int_L
 \sum_{j,k} a_j a_k\theta_j \wedge \varphi_k$$
 and using $\theta_j=\sum_i\lambda_{ij}\alpha_i,
 \varphi_k=\sum_i\mu_{ik}\beta_i$ this is the same as (\ref{fg}).
 \vskip .25cm

 \section{Symplectic aspects}
 
 We have seen that the function $F$ embeds the moduli space of special
 Lagrangian submanifolds of $X$ which are deformations of $L$  as a submanifold
 of  $H^1(L)\times  H^{n-1}(L)$. A vector space  of the form $V\oplus V^*$ also
 has a natural symplectic form $w$ defined by
 $$w((v, \alpha),(v',\alpha'))=\langle v, \alpha'\rangle-\langle v',
 \alpha\rangle$$  so that $H^1(L)\times  H^{n-1}(L)$ may be considered as a
 symplectic manifold. 
 We    shall    now    show    the   following:  \begin{theorem}  The  map  $F$
 embeds $M$ in $H^1(L)\times  H^{n-1}(L)$ as a  Lagrangian submanifold. 
 \end{theorem}  \noindent{\it Proof:} We need to use the algebraic identity
 (iii)  in Section 2 relating
 $\omega$ and $\Omega$ on $X$:
 $$\omega \wedge \Omega=0$$   
 Let $Y$ and $Z$ be two  vector fields, 
then taking interior products with this identity, we obtain

$$0=(\iota(Z)\iota(Y)\omega)\wedge\Omega
-\iota(Y)\omega\wedge\iota(Z)\Omega
+\iota(Z)\omega\wedge\iota(Y)\Omega
+\omega\wedge(\iota(Z)\iota(Y)\Omega)$$
and restricting to a special Lagrangian submanifold $L$, since $\omega$ and
$\Omega$ vanish, we have  
$$\iota(Y)\omega\wedge\iota(Z)\Omega=\iota(Z)\omega\wedge\iota(Y)\Omega$$ 
 Now for $Y$ and $Z$ use vector
fields extending $\partial/\partial t_i$ and  $\partial/\partial t_j$, and we then 
obtain
on $L$ 
 $$\theta_i \wedge \varphi_j=\theta_j\wedge \varphi_i$$
 Thus, integrating,
 $$\int_L \theta_i \wedge \varphi_j=\int_L\theta_j\wedge \varphi_i$$
 and so using $\theta_j=\sum_i\lambda_{ij}\alpha_i,
 \varphi_k=\sum_i\mu_{ik}\beta_i$,
 \begin{equation}
 \sum_i\lambda_{ik}\mu_{ij}=\sum_i\lambda_{ij}\mu_{ik}
 \label{lm}
 \end{equation}
 From the definitions of the coordinates $u$ and $v$ in (\ref{u}) and (\ref{v})
 we  have    $$\lambda_{ij}=\frac{\partial u_i}{\partial t_j},\quad
 \mu_{ij}=\frac{\partial  v_i}{\partial t_j}$$
 so that (\ref{lm}) becomes
 $$ \sum_i\frac{\partial u_i}{\partial t_k}\frac{\partial
 v_i}{\partial t_j}=\sum_i\frac{\partial u_i}{\partial t_j}\frac{\partial
 v_i}{\partial t_k}$$
 But this says precisely that
 $$F^*(\sum_i du_i\wedge dv_i)=0$$
 \vskip .5cm
 It is well-known that a Lagrangian submanifold of the cotangent bundle $T^*N$
 of a manifold for which the projection to $N$ is a local diffeomorphism is
 locally defined as the image of a section $d\phi:N \rightarrow T^*N$ for some
 function  $\phi:N\rightarrow {\bf R}$. Thus, as a consequence of the theorem, 
taking
 $N=H^1(L)$, we can write
 \begin{equation}
 v_j=\frac{\partial \phi}{\partial u_j}
 \label{phi}
 \end{equation}
 for some function $\phi(u_1,\dots,u_m)$. 
 From Proposition 2  the natural metric on $M$ can be written
 in the coordinates $u_1,\dots,u_m$ as 
 \begin{equation}
 g=F^*G=\sum_i du_i dv_i=\sum_{i,j}\frac{\partial^2\phi}{\partial u_i\partial
 u_j}  du_i  du_j  \label{metric}
 \end{equation}

  \vskip .25cm
 Equally, we can take $N=H^{n-1}(L)$ and find a function $\psi(v_1,\dots,v_m)$
 to represent the metric in a similar form:
 $$g=\sum_{ij}\frac{\partial^2\psi}{\partial v_i\partial v_j} dv_i dv_j$$
 The two functions $\phi, \psi$ are related by the classical Legendre transform.
 \vskip .25cm
 \noindent {\it Remark:} Metrics of the above form are said to be of {\it
 Hessian  type}.  V.Ruuska characterized them in \cite{Ru} as those metrics
 admitting an abelian  Lie algebra of gradient vector fields, the local action 
being simply
 transitive.

 \vskip .5cm
 Given that $M$  parametrizes {\it special} Lagrangian submanifolds, it would
 seem reasonable to seek an analogue of the special condition which $M$
 might inherit from the embedding $F$. Now the generators of
 $\Lambda^mV$ and $\Lambda^m V^*$ define two constant $m$-forms $W_1$ and $W_2$
 on  the  $2m$-dimensional manifold $V\times V^*$. We could say that a
 Lagrangian submanifold of $V\times V^*$ is special if a linear combination of
 these forms vanishes, in addition to the symplectic form $w$. With this set-up
 we have:
 \begin{prop} The map $F$ embeds $M$ as a {\it special} Lagrangian submanifold
 if and only if any of the following equivalent statements holds:
 \begin{itemize}
 \item
 $\phi$ satisfies the Monge-Amp\`ere equation $\det(\partial^2\phi/\partial u_i
 \partial u_j)=c$
 \item
 $\psi$ satisfies the Monge-Amp\`ere equation $\det(\partial^2\psi/\partial v_i
 \partial v_j)=c^{-1}$
 \item
 The volume of the torus $H^1(L_t,{\bf R/Z})$ is independent of $t\in M$
 \item
 The volume of the torus $H^{n-1}(L_t,{\bf R/Z})$ is independent of $t\in M$
 \end{itemize}
 \end{prop}
 \noindent{\it Proof:} For the first part, note that, using the coordinates
 $u_1,\dots,u_m$, the $m$-form
 $c_1W_1+c_2W_2$  vanishes on $F(M)$ if and only if
 $$c_1du_1\wedge \dots \wedge du_m+c_2\det(\partial^2\phi/\partial u_i
 \partial u_j)du_1\wedge \dots \wedge du_m=0$$
 which gives
 $$\det(\partial^2\phi/\partial u_i
 \partial u_j)=-c_1/c_2=c$$
 Interchanging the roles of $V$ and $V^*$ gives the second statement.

 To determine the volume of the torus $H^1(L_t,{\bf R/Z})$, we take a
 basis $a_1,\dots,a_m$ of harmonic 1-forms, normalized by
 $$\int_{A_i} a_j=\delta_{ij}$$
 and then the volume is $\sqrt{\det (a_i,a_j)}$ using the inner product on
 harmonic forms. Now from the definition of $\lambda_{ij}$, the normalized
 harmonic forms are
 $$a_j=\sum_k (\lambda^{-1})_{kj}\theta_k$$
 and the inner product
 $$(\theta_j,\theta_k)=\int_L \theta_j\wedge \ast \theta_k=\sum_i
 \lambda_{ij}\mu_{ik}$$
 Thus the volume is 
 $$\sqrt{\det(\mu \lambda^{-1})}$$
 Now in the coordinates $t_1,\dots,t_m$ the form $c_1W_1+c_2W_2$ restricted to
 $F(M)$ is 
 $$(c_1\det \lambda + c_2 \det \mu) dt_1\wedge \dots \wedge dt_m$$
 and this vanishes if and only if $\det(\mu \lambda^{-1})=-c_1/c_2$.
 The final statement follows in similar way. The volume in this case is
 $\sqrt{\det(\lambda    \mu^{-1})}$  \vskip .25cm
 \noindent {\it Remarks:}
 
 \noindent 1.  The relationship between pairs of solutions to the
 Monge-Amp\`ere equations related by the Legendre transform is well-documented
 (see  \cite{Blas}).
 
 \noindent 2. On  any  special  Lagrangian submanifold the volume form is the 
restriction of
 $\Omega_2$, and $\Omega_2$  is  {\it closed} in $X$, so the cohomology class of
 the  volume  form is independent of $t$. Thus the 1-dimensional torus
 $H^n(L,{\bf R/Z})$ has constant volume.
 
 \noindent 3. In the case where $X$ is hyperk\"ahler and $L$ is complex
 Lagrangian with respect to the complex structure $I$, then the flat metric on 
 $H^1(L,{\bf R/Z})$ is K\"ahler and its volume is essentially the Liouville
 volume  of  the  K\"ahler form. But the symplectic form on the torus is
 cohomologically determined: if $[\omega_1]\in H^2(L,{\bf R})$ is the
 cohomology  class  of the $I$-K\"ahler form of $X$, then for $\alpha,\beta
 \in H^1(L,{\bf R})$ the skew form is given by 
 $$\langle \alpha,\beta\rangle [\omega_1]^k=\alpha \wedge \beta \wedge 
 [\omega_1]^{k-1}$$
 Since this is entirely cohomological, it is independent of $t$.
 
 \noindent 4. Another geometrical interpretation of the structure on $M$ is as
 an affine hypersurface $x_{m+1}=\phi(x_1,\dots,x_m)$.  The  Legendre transform  
then corresponds to the dual
 hypersurface of tangent planes, and a solution to the Monge-Amp\`ere equation
 describes a {\it parabolic affine hypersphere} (\cite{Cal2},
 \cite{Blas}).

 \section{ K\"ahler metrics}
 
 The approach of Strominger, Yau and Zaslow takes the moduli space not just of
 special  Lagrangian  submanifolds, but of submanifolds together with 
 flat unitary line bundles (``supersymmetric cycles''). Since a flat line
 bundle on $L$ is classified by an element of $H^1(L, {\bf R/Z})$, then by
 homotopy invariance (we are working locally or on
 a  simply  connected space) this augmented moduli space  can be taken to be 
 $$M^c=M\times H^1(L, {\bf R/Z})$$
 The tangent space $T_m$ at a point of $M^c$ is thus canonically
 $$T_m\cong H^1(L,{\bf R})\oplus H^1(L,{\bf R})\cong H^1(L,{\bf R})\otimes {\bf
 C}$$
This is a complex vector space, so $M^c$ has an almost complex structure.
Moreover, for any real vector space $V$, a positive definite inner product on
$V$ defines a hermitian form on $V\otimes {\bf C}$, so $M^c$ has a hermitian
metric. We then have:
  \begin{prop} The almost complex structure $I$ on $M^c$ is integrable and the
 inner product on $H^1(L,{\bf R})$ defines a K\"ahler metric on $M^c$. 
 \end{prop}
 \noindent{\it Proof:} Use the basis $\alpha_1,\dots,\alpha_m$ of $H^1(L,{\bf
 R})$ to give coordinates $x_1,\dots,x_m$ on the universal covering of the
 torus $H^1(L, {\bf R/Z})$. Then $(t_1,\dots,t_m,x_1,\dots,x_m)$ are local
 coordinates for $M^c$ and from (\ref{coord}) the almost complex structure is
 defined by 
 \begin{eqnarray*}
I( \partial/\partial t_j)&=&\sum_i\lambda_{ij}\partial/\partial x_i\\
I(\sum_i \lambda_{ij}\partial/\partial x_i)&=&-\partial/\partial t_j
\end{eqnarray*}
If we  define the complex vector fields
 $$X_j=\partial/\partial t_j-iI\partial/\partial t_j=\partial/\partial
 t_j-i\sum \lambda_{jk}\partial/\partial x_j$$
 then these satisfy $IX_j=iX_j$ and so form a basis for the $(1,0)$ vector
 fields. The forms $\theta_j$ defined by 
 $$\theta_j= \sum \lambda_{jk}dt_k-idx_j$$
annihilate the $X_j$ and thus form a basis of the $(0,1)$-forms. But from
(\ref{u})
$$\theta_j=d(u_j-ix_j)$$
so that $w_j=u_j+ix_j$ are complex coordinates, and the complex structure is
integrable.
\vskip .25cm
The  2-form $\tilde \omega$ for the Hermitian metric is defined by 
$$\tilde\omega(\partial/\partial t_j,\partial/\partial x_k)=g(\partial/\partial
t_j, I\partial/\partial x_k)$$ and from the definition of $I$, 
$$\tilde\omega(\partial/\partial t_j,\partial/\partial x_k)=-\sum_l
\lambda_{lk}^{-1}g_{jl}$$
But from Proposition 2 the metric is $F^*G$, so in the local coordinates
$t_1,\dots,t_m$,
$$g_{ij}=\sum_k \frac{\partial u_k}{\partial t_i} \frac{\partial v_k}{\partial
t_j}=\sum_k \lambda_{ki}\mu_{kj}$$
(note that symmetry follows from (\ref{lm})). 
Thus,
$$\tilde \omega= -\sum_{j,k} \mu_{kj}dt_j\wedge dx_k=-\sum_k dv_k\wedge dx_k$$
from (\ref{v}). This is clearly closed, so the metric is K\"ahlerian.
\vskip .25cm
\noindent{\it Remark:} Since $v_k=\partial \phi/\partial u_k$, we can also
write
\begin{eqnarray*}
\tilde \omega &=&-\sum_{j,k}(\partial^2 \phi/\partial u_j \partial u_k)
du_j\wedge dx_k\\
&=&(2i)^{-1} \partial \bar \partial \phi
\end{eqnarray*}
so that $\phi/2$ is a K\"ahler potential for this metric. Such metrics, where
the potential depends only on the real part of the complex variables, were 
considered by Calabi in \cite{Cal1}.
\vskip .25cm
We have seen that the pulled-back metric $F^*G$ defines a K\"ahler metric on
$M^c$. If we pull back the constant $m$-form $F^*W_1=du_1\wedge \dots \wedge
du_m$, then this defines directly a complex $m$-form
$$\tilde \Omega^c=d(u_1+ix_1)\wedge \dots d(u_m+ix_m)=dw_1\wedge \dots \wedge
dw_m$$ which is clearly non-vanishing and holomorphic. Using this, we have:
\begin{prop} The holomorphic $m$-form $\tilde \Omega^c$ has constant length
with respect to the K\"ahler metric if and only if any of the equivalent
conditions of Proposition 3 hold.
\end{prop}
\noindent{\it Proof:} First note that
$$dw_j\wedge d\bar w_j=(\sum_k \lambda_{jk}dt_k+idx_j)\wedge(\sum_k
\lambda_{jk}dt_k-idx_j)=2idx_j\wedge\sum_j\lambda_{jk}dt_k$$
Thus 
$$dw_1\wedge\dots dw_m\wedge d\bar w_1\dots \wedge d \bar
w_m=(2i)^m(\det\lambda) dx_1\wedge\dots \wedge dx_m\wedge dt_1\dots \wedge
dt_m$$ But 
$$\tilde \omega^m= (-\sum \mu_{kj}dt_j\wedge dx_k)^m=(-1)^{m(m+1)/2}(\det \mu)
dx_1\wedge\dots \wedge dx_m\wedge dt_1\dots \wedge dt_m$$
Thus $\tilde \Omega^c$ has constant length iff $\det \mu$ is a constant
multiple of $\det \lambda$. But from the proof of Proposition 3, this is
equivalent to the volume of the torus being constant.

\noindent Note that we could equally have argued using the Monge-Amp\`ere
equation for the K\"ahler potential.
\vskip .5cm
\noindent We have thus seen  that if $F$ maps $M$ to a {\it special} Lagrangian
submanifold of $H^1(L)\times H^{n-1}(L)$, the complex manifold $M^c$ has a
natural {\it Calabi-Yau} 
metric.
\vskip .5cm
\noindent{\it Remarks:} 
\vskip .15cm
\noindent 1. It is not hard to see that the tori $H^1(L,{\bf R/Z})\times \{t\}$
in $M^c$ are special Lagrangian with respect to the natural K\"ahler metric and
the holomorphic form $i^m\tilde \Omega ^c$. Since the first Betti number of
this torus 
is $m=\dim M$, the family parametrized by $t\in M$ is complete by McLean's
result, and so we can repeat the process to find another K\"ahler manifold. The 
reader may
easily verify that the roles of $\lambda, \mu$, $u_i,v_i$, $\phi$ and $\psi$
are interchanged. In \cite{SYZ}, one  begins with a Calabi-Yau manifold with a
family of special Lagrangian tori, and produces its ``mirror'' $M^c$ in the
above sense. Performing the process a second time one obtains some sort of
approximation to the first manifold. The metric defined here, however, even
when it is Calabi-Yau, will hardly ever extend to a {\it compact} manifold,
since it has non-trivial Killing fields $\partial/\partial x_i$ -- by Bochner's 
original
Weitzenb\"ock argument, zero Ricci tensor would imply that these are covariant
constant.
 \vskip .15cm \noindent 2. The simplest case of the above process consists of 
considering
elliptic curves in a hyperk\"ahler 4-manifold (a 2-dimensional Calabi-Yau
manifold). Thus $m=2$ and we obtain a 4-dimensional hyperk\"ahler metric on
$M^c$. The existence of two Killing fields shows that it must be produced from
the Gibbons-Hawking ansatz \cite{GH} using a harmonic function of two
variables. From the above arguments, this means that the 2-dimensional
Monge-Amp\`ere equation can be reduced to  Laplace's equation in two
variables. In fact, as the reader will find in \cite{Dar}, this is classically
 known. In the same way curves of  genus $g$ in (for example) a K3 surface
 generate  a  solution to the $2g$-dimensional Monge-Amp\`ere equation.

\end{document}